\documentclass{amsart}
\usepackage{amsmath, amssymb}

\newtheorem{thm}{Theorem}

\theoremstyle{remark}
\newtheorem{rmk}[thm]{Remark}

\numberwithin{thm}{section}
\numberwithin{equation}{section}

\newcommand{\pv}{\mathbf{P}}
\newcommand{\cv}{\mathbf{C}}
\newcommand{\nv}{\mathbf{N}}
\newcommand{\qv}{\mathbf{Q}}
\newcommand{\ord}{\textup{ord}}
\def\Re{{\sf Re}\,}

\begin{document}

\title{A brief survey on local holomorphic dynamics in higher dimensions}

\author{Feng Rong}

\address{Department of Mathematics, Shanghai Jiao Tong University, 800 Dong Chuan Road, Shanghai, 200240, P.R. China}
\email{frong@sjtu.edu.cn}

\thanks{The author is partially supported by the National Natural Science Foundation of China (Grant No. 11371246).}



\begin{abstract}
We give a brief survey on local holomorphic dynamics in higher dimensions. The main novelty of this note is that we will organize the material by the ``level" of local invariants rather than the type of maps.
\end{abstract}

\maketitle

\section{Introduction}

Let $f$ be a holomorphic map in $\cv^n$ with a fixed point, which we assume to be the origin. The local (discrete) holomorphic dynamics studies the asymptotic behavior of $f$ in a neighborhood of the fixed point under iterations. There are several well-written surveys on this subject, see e.g. \cite{A:Survey, A:Discrete, B:Local}. The aim of this short note is twofold. First, we will present some more recent results in this area which were not covered in previous surveys. Second, we organize the material in a different way than before so as to emphasize the importance of the ``third-level" local invariants in future studies.

A quantity will be called a \textit{local invariant} if it only depends on the map $f$, i.e. invariant under changes of local coordinates. (\'Ecalle \cite{Ec} gave a detailed study on local invariants of holomorphic maps, although the dynamics associated with these invariants are not clear.) We will divide local invariants into three levels, depending on from which ``level" of the Taylor expansion of $f$ at $0$ the invariants are defined. Roughly speaking, a \textit{first-level} invariant comes from the linear part of the Taylor expansion; a \textit{second-level} invariant comes from the leading nonlinear part of the Taylor expansion; and a \textit{third-level} invariant comes from higher order nonlinear part of the Taylor expansion. When defining these local invariants, we will always use some suitable local coordinates. However, all these local invariants have been shown to be well-defined, i.e. independent of the choice of (allowable) local coordinates.

There are two typical types of results in local holomorphic dynamics. One type is to give normal forms or even linearizations via conjugations, the other is to describe the attracting set of a given map.

Two maps $f$ and $g$ are said to be (holomorphically) \textit{conjugate} if there exists a biholomorphism $\varphi$ such that $f\circ\varphi=\varphi\circ g$. Obviously, if $f$ and $g$ are conjugate then their local dynamics are equivalent. For a given map $f$, the ``simplest" $g$ it is conjugate to is called its \textit{normal form}. The best one can hope for is that $g$ is the linear part of $f$, in which case we say that $f$ is \textit{linearizable}. The well-known Poincar\'{e}-Dulac normal form and Brjuno's linearization theorem are typical examples. This type of results are certainly important. However, the majority of the note will be devoted to results on the attracting sets.

A point $p$ in a neighborhood of $0$ is in the \textit{attracting set} of $f$ if $f^k(p)$ goes to $0$ as $k$ goes to the infinity. Here, of course, $f^k$ stands for the $k$-th iteration of $f$. The ultimate goal of the local dynamical study is to give a complete description of the asymptotic behavior of a map in a \textit{full} neighborhood of the fixed point. The first step in achieving this goal is to give a complete study on the attracting set. The well-known Leau-Fatou Flower Theorem is the ``model" result. Much of the work on the attracting sets in higher dimensions can be viewed as generalizations of the Leau-Fatou Flower Theorem.

Due to the limit of space, the results surveyed in the note are by no means complete. Our focus will be on results obtained in the past few years. For more detailed information on earlier results and results in one dimension, please see the existing surveys cited above. We would like to thank Filippo Bracci and the referee for valuable comments.

\section{The first-level invariants}

\subsection{The multipliers}

Let $f$ be a holomorphic map in $\cv^n$ with the origin as a fixed point. Write $f$ as
$$f(z)=L(z)+P_2(z)+P_3(z)+\cdots,$$
where $L(z)$ is the linear part of $f(z)$ and $P_k(z)$ are homogeneous of degree $k$, $k\ge 2$. Write $L(z)=L\cdot z$, where $L$ is an $n\times n$ matrix. The \textit{multipliers} of $f$ are defined to be the eigenvalues $\{\lambda_j\}_{j=1}^n$ of $L$.

Denote by $\nv$ the set of non-negative integers. For $\alpha=(\alpha_1,\cdots,\alpha_n)\in \nv^n$, set $\lambda^\alpha=\prod_{j=1}^n \lambda_j^{\alpha_j}$ and $|\alpha|=\sum_{j=1}^n \alpha_j$. Define
$$\omega(m)=\min\limits_{2\le|\alpha|\le m}\min\limits_{1\le j\le n}|\lambda^\alpha-\lambda_j|,\qquad m\ge 2.$$
We say that the multipliers of $f$ satisfy the \textit{Brjuno condition} if
$$\sum_{i\ge 0} \frac{1}{p_i}\log \frac{1}{\omega(p_{i+1})}<\infty,$$
where $\{p_i\}_{i=0}^{\infty}$ is a sequence of integers with $1=p_0<p_1<\cdots$. The following is the best known linearization result (improving earlier results by Siegel \cite{S:1942}).

\begin{thm}[Brjuno, \cite{Br:1971-72}]\label{T:Brjuno}
Let $f$ be a holomorphic map in $\cv^n$ with the origin as a fixed point. If $df_0$ is diagonalizable and the multipliers of $f$ satisfy the Brjuno condition, then $f$ is holomorphically linearizable.
\end{thm}

A \textit{resonance} for $f$ is a relation of the form
$$\lambda^\alpha-\lambda_j=0,\qquad |\alpha|\ge 2,\ 1\le j\le n.$$
Obviously, Theorem \ref{T:Brjuno} is not applicable in the presence of resonances. A map $f$ is said to be \textit{quasi-parabolic}, if $df_0$ is diagonalizable and $\lambda_j=1$ for $1\le j\le m<n$ and $\lambda_j\neq1$ but $|\lambda_j|=1$ for $m+1\le j\le n$. Note that quasi-parabolic maps always have resonances. However, inspired by a partial linearization result by P\"{o}schel \cite{Po:Manifold}, the author proved the following linearization result for quasi-parabolic maps.

\begin{thm}[Rong, \cite{R:Quasi1}]
Let $f$ be a holomorphic map in $\cv^n$ with the origin as a quasi-parabolic fixed point. Assume that there exists an $m$-dimensional manifold $M$ of fixed points through $0$ such that $df_p=df_0$ for every $p\in M$. If $\{\lambda_j\}_{j=m+1}^n$ satisfy the Brjuno condition, then $f$ is holomorphically linearizable.
\end{thm}

This was later generalized to more general settings by Raissy \cite{Ra:Brjuno}.

\subsection{Multi-resonance}

As we have seen above, the presence of resonances is usually an obstacle for the local dynamical study. However, in a recent development, the presence of resonances has been used in a positive way to study the attracting sets of certain maps.

Assume that $df_0$ is diagonalizable and that there are resonances among the multipliers $\{\lambda_j\}_{j=1}^n$. If the resonances are generated over $\nv$ by a finite number of $\qv$-linearly independent multi-indices, $f$ is said to be \textit{multi-resonant}. In \cite{BZ:One}, Bracci and Zaitsev studied \textit{one-resonant} maps and obtained sufficient conditions for the existence of local attracting basins. This was later generalized to multi-resonant maps by Bracci, Raissy and Zaitsev \cite{BRZ:Multi}. More recently, Raissy and Vivas \cite{RaV} gave a more detailed study on \textit{two-resonant} maps, and Bracci and the author \cite{BR:Quasi} studied \textit{quasi-parabolic one-resonant} maps.

The basic idea of this line of study is as follows: first by using the multi-resonance, $f$ can be projected into a lower-dimensional map $\hat{f}$, the so-called \textit{parabolic shadow} of $f$, such that $\hat{f}$ is \textit{tangent to the identity} at the origin; then using the local attracting basin of $\hat{f}$ and some attracting conditions on the ``fibers," one can create an attracting basin for $f$. Recall that a holomorphic map $f$ is said to be tangent to the identity at $0$ if $df_0=I_n$, the identity matrix.

\subsection{Diagonalization}

Most of the results in local holomorphic dynamics assume the linear part of the maps under study to be diagonalizable. However, the non-diagonalizable case is certainly important and worth studying. For instance, Yoccoz \cite{Y:Brjuno} pointed out that the Brjuno condition is in general not sufficient for holomorphic linearization in the non-diagonalizable case (see also \cite{DG:Jordan}).

The method of \textit{blow-up} is a very important tool in the study of local holomorphic dynamics. It is particularly so for the study on attracting sets in the non-diagonalizable case. Indeed, Abate \cite{A:Diag} gave an explicit description of how to systematically diagonalize a non-diagonalizable map via blow-ups.

There are very few results in the non-diagonalizable case, see e.g. \cite{A:Jordan}. Recently, the author \cite{R:Jordan} gave a somewhat systematic study of the non-diagonalizable case in dimension two. In particular, the attracting basin studied in \cite{A:Jordan} was recovered as a special case.

\section{The second-level invariants}

\subsection{The order and characteristic directions}

Let us first recall the well-known \textit{Leau-Fatou Flower Theorem} from the one-dimensional theory (see e.g. \cite{M:One}).

\begin{thm}[Leau, \cite{Leau}; Fatou, \cite{Fatou}]
Let $f$ be a holomorphic map in $\cv$ with the origin as a fixed point. Assume that $f$ is tangent to the identity with order $\nu$, i.e. $f$ can be written as
$$f(z)=z+az^\nu+O(z^{\nu+1}),\qquad \nu\ge 2,\ a\neq 0.$$
Then there exist $\nu-1$ ``attracting petals" for $f$ at the origin.
\end{thm}

A central theme in the study on attracting sets for holomorphic maps in higher dimensions has been to generalize the Leau-Fatou Flower Theorem. To state the known results so far, let us first make several definitions.

Let $f$ be a holomorphic map in $\cv^n$, tangent to the identity at the origin. Write $f$ as
$$f(z)=z+P_2(z)+P_3(z)+\cdots,$$
where $P_k(z)$, $k\ge 2$, are $n$-tuples of homogeneous polynomials of degree $k$. The \textit{order} $\nu$ of $f$ is defined as
$$\nu:=\min\{k:\ P_k(z)\not\equiv 0\}.$$
Write $z=(z_1,\cdots,z_n)$ and $P_\nu(z)=(P_{\nu,1}(z),\cdots,P_{\nu,n}(z))$ and denote by $[\cdot]$ the canonical projection from $\cv^n\backslash\{0\}$ to $\pv^{n-1}$. A direction $[v]=[z_1:\cdots:z_n]$ is called a \textit{characteristic direction} of $f$ if
$$(P_{\nu,1}(z),\cdots,P_{\nu,n}(z))=\lambda(z_1,\cdots,z_n),\qquad \lambda\in\cv.$$
If $\lambda\neq 0$, then $[v]$ is said to be \textit{non-degenerate}, otherwise \textit{degenerate}.

An \textit{attracting petal} of dimension $d$ for $f$ at the origin is an injective holomorphic map $\varphi:\Delta\rightarrow \cv^n$ satisfying the following properties:\\
(i) $\Delta$ is a simply connected domain in $\cv^d$ with $0\in \partial\Delta$;\\
(ii) $\varphi$ is continuous on $\partial\Delta$ and $\varphi(0)=0$;\\
(iii) $\varphi(\Delta)$ is invariant under $f$ and $f^k(\varphi(\zeta))\rightarrow 0$ as $k\rightarrow \infty$ for any $\zeta\in \Delta$.\\
Furthermore, if $[\varphi(\zeta)]\rightarrow [v]\in \pv^{n-1}$ as $\zeta\rightarrow 0$, then $\varphi$ is said to be \textit{tangent to} $[v]$ at $0$. If there are $\nu-1$ attracting petals tangent to $[v]$ at $0$, then we say they form an \textit{attracting flower} tangent to $[v]$ at $0$. When $d=1$, an attracting petal is also called a \textit{parabolic curve}. When $1<d<n$, an attracting petal is also called a \textit{parabolic manifold}. When $d=n$, an attracting petal is a (parabolic) \textit{attracting domain}.

The first main generalization of the Leau-Fatou Flower Theorem to higher dimensions is the following

\begin{thm}[\'Ecalle, \cite{Ec}; Hakim, \cite{H:Parabolic}]\label{T:Hakim}
Let $f$ be a holomorphic map in $\cv^n$, tangent to the identity at the origin. If $f$ is of order $\nu<\infty$ and $[v]$ is a non-degenerate characteristic direction of $f$, then there exists a one-dimensional attracting flower of $f$ tangent to $[v]$ at $0$.
\end{thm}

A similar result holds for quasi-parabolic maps, which was proven in dimension two by Bracci and Molino \cite{BM:Quasi} and in higher dimensions by the author \cite{R:Quasi2}. To be more precise, we need some definitions.

Let $f$ be a quasi-parabolic map in $\cv^n$ with eigenvalue 1 of multiplicity $l$ and other eigenvalues $\lambda_j$, $1\le j\le m=n-l$. Set $\Lambda=\textup{Diag}(\lambda_1,\cdots,\lambda_m)$. In a suitable local coordinates $(z,w)\in \cv^l\times \cv^m$, we can then write $f$ as
$$\left\{\begin{aligned}
z_1&=z+p(z)+r(z,w),\\
w_1&=\Lambda w+q(z)+s(z,w),
\end{aligned}\right.$$
where $p,q,r,s$ are all of degree greater or equal to two and $r(z,0)=s(z,0)=0$.

We say that $f$ is in \textit{ultra-resonant} form if $\ord p(z)\le \ord q(z)$, in which case we call $\nu=\ord p(z)$ the \textit{order} of $f$. Assume that $\nu<\infty$, and let $p_\nu(z)$ be the lowest order term of $p(z)$. A \textit{characteristic direction} of $f$ is of the form $[v]=[z_1:\cdots:z_l:0:\cdots:0]$ where $[u]=[z_1:\cdots:z_l]$ is a characteristic direction of $p_\nu(z)$, i.e. $p_\nu(z)=\lambda z$ for some $\lambda\in \cv$. And $[v]$ is said to be \textit{non-degenerate} if $\lambda\neq 0$, otherwise \textit{degenerate}.

If $f$ has a characteristic direction $[v]$, then in a suitable local coordinates it can be assumed that $[v]=[1:0:\cdots:0]$. Write $z=(x,y)\in \cv\times\cv^{l-1}$. Set $\mu=\min\{k;\ x^kw_j \textup{ in } s(z,w)\}$. We say that $f$ is \textit{dynamically separating} in $[v]$ if $\mu\ge \nu-1$.

\begin{thm}[Bracci-Molino, \cite{BM:Quasi}; Rong, \cite{R:Quasi2}]
Let $f$ be a holomorphic map in $\cv^n$, with a quasi-parabolic fixed point at the origin. If $f$ is of order $\nu<\infty$, has a non-degenerate characteristic direction $[v]$, and $f$ is dynamically separating in $[v]$, then there exists a one-dimensional attracting flower of $f$ tangent to $[v]$ at $0$.
\end{thm}

\subsection{The director and residual index}

Let $f$ be a holomorphic map in $\cv^n$, tangent to the identity at the origin. Assume that $f$ has order $\nu<\infty$ and has a non-degenerate characteristic direction $[v]$. In suitable local coordinates $z=(x,y)\in\cv\times\cv^{n-1}$, it can be assumed that $[v]=[1:0]$. Write $f$ as
\begin{equation}\label{E:f}
\left\{\begin{aligned}
x_1&=x+p_\nu(x,y)+O(\nu+1),\\
y_1&=y+q_\nu(x,y)+O(\nu+1),
\end{aligned}\right.
\end{equation}
where $p_\nu(x,y)$ and $q_\nu(x,y)$ are homogeneous of degree $\nu$.

Under the blow-up $y=xu$, the blow-up map $\tilde{f}$ takes the form
\begin{equation}\label{E:f1}
\left\{\begin{aligned}
x_1&=x+x^\nu p_\nu(1,u)+O(x^{\nu+1}),\\
u_1&=u+x^{\nu-1} r(u)+O(x^\nu),\\
\end{aligned}\right.
\end{equation}
where $r(u)=q_\nu(1,u)-p_\nu(1,u)u$. The matrix
$$A:=p_\nu^{-1}(1,0)r'(0)$$
is a local invariant associated with $f$ and its $n-1$ eigenvalues are called the \textit{directors} of $f$ in the non-degenerate characteristic direction $[v]$.

\begin{thm}[Hakim, \cite{H:Parabolic1}]\label{T:Hakim1}
Let $f$ be a holomorphic map in $\cv^n$, tangent to the identity at the origin. Assume that $f$ has order $\nu<\infty$ and has a non-degenerate characteristic direction $[v]$. Let $\alpha_i$, $1\le i\le n-1$, be the directors of $f$ in $[v]$. Suppose that there exists $\alpha>0$ such that $\Re\alpha_j>\alpha$ for $1\le j\le l$ and $\Re\alpha_{l+k}<\alpha$ for $1\le k\le n-1-l$. Then there exists an $(l+1)$-dimensional attracting flower of $f$ tangent to $[v]$ at $0$.
\end{thm}

Similar results hold for quasi-parabolic maps (see \cite{R:Quasi3}) and semi-attractive maps (see e.g. \cite{Fatou1, U:Local, H:Semi, Ri:Semi, R:Semi}). Recall that a holomorphic map $f$ is said to be \textit{semi-attractive} at $0$ if $df_0=\textup{Diag}(I_l,A)$, where $I_l$ is the identity matrix with size $1\le l\le n-1$ and the eigenvalues of $A$ all have modulus less than one.

Theorems \ref{T:Hakim} and \ref{T:Hakim1} deal with holomorphic maps tangent to the identity with a non-degenerate characteristic direction, which is a generic condition. It would be desirable to obtain a ``full" generalization of the Leau-Fatou Flower Theorem without this generic condition. So far, this has only been achieved in dimension two by Abate \cite{A:Residual}.

\begin{thm}[Abate, \cite{A:Residual}]\label{T:Abate}
Let $f$ be a holomorphic map in $\cv^n$. Assume that the origin is an isolated fixed point of $f$ and $f$ is tangent to the identity at $0$. Then there exists a one-dimensional attracting flower of $f$ at $0$.
\end{thm}

The main point of the proof of Theorem \ref{T:Abate} is to show that after a sequence of blow-ups one gets a blow-up map which is generic, i.e. with a non-degenerate characteristic direction, as in Theorem \ref{T:Hakim}. For this purpose, Abate introduced a key local invariant, the \textit{residual index}, defined as follows.

Let $\tilde{f}$ be a holomorphic map in $\cv^2$, tangent to the identity at the origin. Assume that there is a line $S$ of fixed points of $\tilde{f}$ through $0$. In local coordinates $(x,u)\in \cv\times\cv$, such that $S$ is given by $\{x=0\}$, we can write $\tilde{f}$ as
$$\left\{\begin{aligned}
x_1&=x+x^\nu p(x,u),\\
u_1&=u+x^\mu q(x,u),\\
\end{aligned}\right.$$
where $p(0,u)\not\equiv 0$ and $q(0,u)\not\equiv 0$, $\nu\ge 2$ and $\mu\ge 1$.

Define a meromorphic function, the \textit{residual function}, $\kappa(u)$ by
$$\kappa(u):=\lim\limits_{x\rightarrow 0} \frac{p(x,u)}{x^{\mu-\nu+1}q(x,u)}.$$
If $\mu<\nu-1$, then $\kappa(u)\equiv 0$. If $\mu>\nu-1$, then $\kappa(u)\equiv \infty$. If $\mu=\nu-1$, then $\kappa(u)=p(0,u)/q(0,u)$. The \textit{residual index} $\iota_0(\tilde{f},S)$ of $\tilde{f}$ at $0$ along $S$ is defined as $\textup{Res}(\kappa(u);0)$. Note that if $\tilde{f}$ is the blow-up map of a holomorphic map tangent to the identity in a non-degenerate characteristic direction and $S$ is the exceptional divisor, then the residual index is exactly the reciprocal of the director as defined above.

Although Theorem \ref{T:Abate} gives a Leau-Fatou Flower Theorem in dimension two, it leaves open two questions: 1. What happens if the origin is non-isolated? 2. Given a \textit{degenerate} characteristic direction, is there always an attracting petal tangent to it?

For results related to the first question, see e.g. \cite{B:Fixed, De:Fixed}. Note that the residual index theorems used in \cite{A:Residual, B:Fixed} have been developed systematically by Abate, Bracci and Tovena \cite{ABT} to much more general settings and also to higher dimensions. It would be desirable to find an effective use of such more general index theorems to the study of local holomorphic dynamics.

For the second question, Abate \cite{A:Residual} already showed that the answer is yes if the residual index of the blow-up map at the given direction along the exceptional divisor does not belong to $\qv^+$. This result was later generalized by Molino \cite{Mo:Index} to the case of non-vanishing residual index (under a mild ``regularity" condition).

In dimension two, using the residual function $\kappa(u)$ defined above, the characteristic directions can be divided into \textit{three} types (cf. \cite{AT:Poincare}).

Let $f$ be a holomorphic map in $\cv^2$, tangent to the identity at the origin. Assume that $f$ has a characteristic direction, which we assume to be $[1:0]$. Then we can write $f$ as in \eqref{E:f}, and its blow-up map $\tilde{f}$ as in \eqref{E:f1}. If $r(u)\equiv 0$, then $f$ is said to be \textit{dicritical} at $0$.

Assume that $f$ is not dicritical at $0$. Then the residual function is given by $\kappa(u)=p_\nu(1,u)/r(u)$. If $0$ is a simple pole of $\kappa(u)$, then $[1:0]$ is a \textit{Fuchsian} characteristic direction of $f$. If $0$ is a pole of $\kappa(u)$ of order greater than one, then $[1:0]$ is an \textit{irregular} characteristic direction of $f$. If $\kappa(u)\equiv 0$ or if $0$ is a removable singularity of $\kappa(u)$, then $[1:0]$ is an \textit{apparent} characteristic direction of $f$.

\begin{thm}[Vivas, \cite{V:Degenerate}]
Let $f$ be a holomorphic map in $\cv^2$, tangent to the identity at the origin. Assume that $f$ has an irregular characteristic direction $[v]$. Then there exists an attracting domain of $f$ tangent to $[v]$ at $0$.
\end{thm}

Vivas \cite{V:Degenerate} also gave sufficient conditions (in terms of the residual index) for the existence of attracting domains in Fuchsian characteristic directions, and studied examples of apparent characteristic directions. See also \cite{V:F-B, Lapan} for related results.

\subsection{The non-dicritical order}

Let $f$ be a holomorphic map in $\cv^n$, tangent to the identity at the origin. Assume that $f$ has order $\nu<\infty$ and has a non-degenerate characteristic direction $[v]$. Let $\alpha_i$, $1\le i\le n-1$, be the directors of $f$ in $[v]$. If $\Re\alpha_j>0$ for $1\le j\le l$ and $\Re\alpha_{l+k}<0$ for $1\le k\le m=n-1-l$, then by Theorem \ref{T:Hakim1} there exists an $(l+1)$-dimensional attracting flower of $f$ tangent to $[v]$ at $0$. In fact, from \cite{H:Parabolic1, AR:Hakim}, it follows that $l+1$ is the maximal dimension of attracting flowers in this case. It is then natural to ask what happens when $\Re\alpha_{l+k}=0$ for all $1\le k\le m$.

In suitable local coordinates $(x,y,z)\in \cv\times\cv^l\times\cv^m$, we can assume that $[v]=[1:0:0]$. Under the blow-up ($y=xu$, $z=xv$), the blow-up map $\tilde{f}$ can be written as (after possible scaling and suitable linear transformations)
\begin{equation}\label{E:f2}
\left\{\begin{aligned}
x_1&=(1-x^{\nu-1})x+O(x^\nu \|w\|,x^{\nu+1}),\\
u_1&=(I_l-x^{\nu-1}B)u+O(x^{\nu-1}\|w\|^2,x^\nu),\\
v_1&=(I_m-x^{\nu-1}C)v+O(x^{\nu-1}\|w\|^2,x^\nu),
\end{aligned}\right.
\end{equation}
where $w=(u,v)$, $B$ is an $l\times l$ matrix with eigenvalues $\alpha_j$, $1\le j\le l$, and $C$ is an $m\times m$ matrix with eigenvalues $\alpha_{l+k}$, $1\le k\le m$.

Rewrite $v_1$ in \eqref{E:f2} as
$$v_1=v+x^{\nu-1}\sum_{|k|=1}^{\nu+1} v^k \gamma_k +O(x^{\nu-1}\|u\|\|w\|,x^\nu),$$
where $k=(k_1,\cdots,k_m)\in \nv^m$ is a multi-index, $|k|=k_1+\cdots+k_m$, and $\gamma_k\in \cv^m$. Then the \textit{non-dicritical order} of $f$ in the characteristic direction $[v]$ is defined as
$$\tau:=\min \{|k|-1;\ \gamma_k\neq 0\}.$$

The name ``non-dicritical" refers to the fact that if $l=0$ then $f$ is dicritical at the origin if and only if all $\gamma_k$ vanish. It will always be assumed that some $\gamma_k\neq 0$, in which case $0\le \tau\le \nu$ (see e.g. \cite{Brochero} for a study in the dicritical case).

\begin{thm}[Rong, \cite{R:Non-dicritical}]
Let $f$ be a holomorphic map in $\cv^n$, tangent to the identity at the origin, with a non-degenerate characteristic direction $[v]$. Let $\{\alpha_j\}_{j=1}^{n-1}$ be the directors of $f$ in $[v]$ and assume that $\Re\alpha_j\le0$ for some $j$. Let $\tau$ be the non-dicritical order of $f$ in $[v]$. If $\tau\ge 1$, then there exists an attracting domain of $f$ tangent to $[v]$ at $0$.
\end{thm}

When $\tau=0$, it is possible for $f$ to admit a ``spiral domain" at the origin (see \cite{R:Non-dicritical} for more details). Note also that in dimension two, if $\tau\ge 1$ then $[v]$ is an irregular characteristic direction of $f$.

\section{The third-level invariants}

\subsection{Essentially non-degenerate}

From the discussion above, it is clear that one of the main problems in the study of local holomorphic dynamics of maps tangent to the identity is to understand the dynamics in degenerate characteristic directions. So far there are only very few results, and only in dimension two (see e.g. \cite{R:Degenerate, V:Degenerate}).

Let $f$ be a holomorphic map in $\cv^2$, tangent to the identity at the origin. Assume that $f$ has order $\nu<\infty$ and $[v]=[x:y]=[1:0]$ is a degenerate characteristic direction of $f$. Write $f$ as
$$\left\{\begin{aligned}
x_1&=x+yp_{\nu-1}(x,y)+O(\nu+1),\\
y_1&=y+yq_{\nu-1}(x,y)+O(\nu+1),
\end{aligned}\right.$$
where $p_{\nu-1}(x,y)$ and $q_{\nu-1}(x,y)$ are homogeneous of degree $\nu-1$. We say that $[1:0]$ is a \textit{generically degenerate} characteristic direction of $f$ if $q_{\nu-1}(x,0)\not\equiv 0$. Note that a generically degenerate characteristic direction is an apparent characteristic direction.

Let $G$ be the group of local changes of coordinates which preserves $[1:0]$ as the degenerate characteristic direction. For each $\phi\in G$, write $(x_\phi,y_\phi)=\phi(x,y)$. Then $f$ is transformed under $\phi$ as
$$\left\{\begin{aligned}
x_{\phi,1}&=x_\phi+P_\phi(x_\phi)+y_\phi O(\nu-1),\\
y_{\phi,1}&=y_\phi+Q_\phi(x_\phi)+y_\phi R_\phi(x_\phi)+y_\phi^2 O(\nu-2),
\end{aligned}\right.$$
where $\ord P_\phi\ge \nu+1$, $\ord Q_\phi\ge \nu+1$ and $\ord R_\phi=\nu-1$. The \textit{essential order} of $f$ in $[1:0]$ is defined as
$$\mu:=\max_{\phi\in G} \ord P_\phi(x_\phi).$$
We say that $[1:0]$ is an \textit{essentially non-degenerate} characteristic direction of $f$ if $\mu<\infty$.

\begin{rmk}
The above definition is slightly different than the original definition given in \cite{R:Degenerate}, where $\mu$ is required to be less than the so-called \textit{virtual order}. However, using simple linear transformations of the form ($X=x$, $Y=y+\alpha x^k$) for $k\ge 2$, it is easy to check that the virtual order is always greater than $\mu$ if $\mu<\infty$.
\end{rmk}

Now assume that $\mu<\infty$ and rewrite $f$ as
$$\left\{\begin{aligned}
x_1&=x-ax^\mu+P(x)+yO(\nu-1),\\
y_1&=y-byx^{\nu-1}+Q(x)+yR(x)+y^2O(\nu-2),
\end{aligned}\right.$$
where $a,b\neq 0$, $\ord Q>\ord P\ge \mu+1$ and $\ord R\ge \nu$. Then the \textit{director} of $f$ in $[1:0]$ is defined as
$$\alpha:=\frac{b}{a^{(\nu-1)/(\mu-1)}}.$$

\begin{thm}[Rong, \cite{R:Degenerate}]
Let $f$ be a holomorphic map in $\cv^2$, tangent to the identity at the origin. Assume that $f$ has an essentially non-degenerate characteristic direction $[v]$ and the director of $f$ in $[v]$ is $\alpha$. If $\Re \alpha>0$, then there exists an attracting domain of $f$ tangent to $[v]$ at $0$.
\end{thm}

Note that the above theorem gives an attracting petal instead of an attracting flower. As first observed by the author \cite{R:Quasi2}, this ``symmetry break-down" is indeed expected in degenerate characteristic directions.

\subsection{Non-dynamically-separating}

Let $f$ be a holomorphic map in $\cv^2$, quasi-parabolic at the origin. Then $[1:0]$ is the only characteristic direction of $f$. Write $f$ as
$$\left\{\begin{aligned}
z_1&=z+P(z)+wS(z,w),\\
w_1&=\lambda w+Q(z)+wR(z)+w^2T(z,w),
\end{aligned}\right.$$
with $|\lambda|=1$ and $\lambda\neq 1$. Set $p=\ord P$, $q=\ord Q$ and $r=\ord R$. Assume that $f$ is in ultra-resonant form, i.e. $q\ge p$, and non-dynamically-separating, i.e. $r<p-1$.

Now rewrite $f$ as
$$\left\{\begin{aligned}
z_1&=z-az^p+O(z^{p+1},w),\\
w_1&=\lambda w-bwz^r+O(z^p,wz^{r+1},w^2),
\end{aligned}\right.$$
with $a,b\neq 0$. We call $p$ the \textit{essential order} and $r$ the \textit{generic order}. The \textit{director} of $f$ is then defined as
$$\alpha:=\frac{b}{\lambda a^{r/(p-1)}}.$$

\begin{thm}[Rong, \cite{R:Quasi5}]
Let $f$ be a holomorphic map in $\cv^2$, quasi-parabolic at the origin. Assume that $f$ is non-dynamically-separating and the director of $f$ is $\alpha$. If $\Re \alpha>0$, then there exists an attracting domain of $f$ tangent to $[1:0]$ at $0$.
\end{thm}

\end{document}